\documentclass[conference]{IEEEtran}
\IEEEoverridecommandlockouts
\usepackage{amsmath,amssymb,amsfonts,mathtools, mathrsfs}
\usepackage{dblfloatfix}
\usepackage{graphicx}
\usepackage{textcomp}
\usepackage[x11names]{xcolor}
\usepackage{xcolor}
\usepackage{todonotes}
\usepackage{tikz}
\usetikzlibrary{shapes.geometric, arrows, positioning, fit, arrows.meta, shapes.misc}
\usepackage{pgfplots}
\pgfplotsset{compat=1.18}
\usepackage{svg}
\usepackage[caption=false]{subfig}
\usepackage{accents}
\usepackage{siunitx}
\DeclareSIUnit{\euro}{\text{€}} 
\usepackage{multirow}
\usepackage{algorithm}
\usepackage{algpseudocode}
\usepackage{pdfpages}
\usepackage{booktabs} 
\usepackage{blkarray}
\usepackage{bigstrut}
\usepackage{float}
\usepackage{makecell}
\usepackage{nicematrix}
\usepackage{xfp}

\definecolor{IEE_light_blue}{HTML}{1E90FF}
\definecolor{IEE_blue}{HTML}{206173}
\definecolor{IEE_red}{HTML}{F70146}
\definecolor{IEE_green}{HTML}{78BE73}
\definecolor{IEE_orange}{HTML}{D58E00}
\definecolor{red1}{HTML}{fdccda}
\definecolor{red2}{HTML}{f70146}

\definecolor{oiOrange}{HTML}{E69F00}
\definecolor{oiSkyBlue}{HTML}{56B4E9}
\definecolor{oiBluishGreen}{HTML}{009E73}
\definecolor{oiYellow}{HTML}{F0E442}
\definecolor{oiBlue}{HTML}{0072B2}
\definecolor{oiVermilion}{HTML}{D55E00}
\definecolor{oiReddishPurple}{HTML}{CC79A7}
\definecolor{oiBlack}{HTML}{000000}

\definecolor{NPAP_220}{HTML}{d4ac82}
\definecolor{NPAP_380}{HTML}{5cc4a4}
\definecolor{NPAP_transformer}{HTML}{f0e874}

\definecolor{NPAP_220}{HTML}{D4AC82}
\definecolor{NPAP_380}{HTML}{5CC4A4}
\definecolor{NPAP_transformer}{HTML}{F0E874}
\definecolor{NPAP_220_inv}{HTML}{A07550}
\definecolor{NPAP_380_inv}{HTML}{02845F}
\definecolor{NPAP_transformer_inv}{HTML}{C8BE20}
\definecolor{NPAP_uniform}{HTML}{6A6A6A}

\newcommand{\databarLine}[2][oiBlue]{%
    \begin{tikzpicture}[baseline=(val.base)] 
        \useasboundingbox (-1.12, -0.05) rectangle (0, 0.3);
        \fill[#1!40] (0, -0.05) rectangle (-#2*0.0004, 0.3);
        \node[anchor=west, inner sep=0pt, text=black] (val) at (-1.12, 0.1) {#2};
    \end{tikzpicture}%
}

\newcommand{\databarTrafo}[2][oiBlue]{%
    \begin{tikzpicture}[baseline=(val.base)] 
        \useasboundingbox (-1.12, -0.05) rectangle (0, 0.3);
        \fill[#1!40] (0, -0.05) rectangle (-#2*0.082, 0.3); 
        \node[anchor=west, inner sep=0pt, text=black] (val) at (-1.12, 0.1) {#2};
    \end{tikzpicture}%
}

\usepackage{hyperref}
\usepackage{csquotes}
\usepackage[utf8]{inputenc}

\usepackage[
style=ieee, 
citestyle=numeric-comp,
backend=biber, 
giveninits=true, 
maxnames=2,
sortcites=true,
]{biblatex}
\AtEveryBibitem{%
    \clearfield{month}%
    \clearfield{language}
    \clearfield{urlyear}%
    \clearfield{issn}%
    \clearfield{doi}
}

\defbibenvironment{bibliography}
  {\list
     {\printfield[labelnumberwidth]{labelnumber}}
     {%
       \setlength{\labelwidth}{2em}
       \setlength{\labelsep}{0.7em}  
       \setlength{\leftmargin}{\dimexpr\labelwidth+\labelsep\relax}
       \setlength{\itemsep}{0.1\baselineskip}
       \setlength{\parsep}{0pt}
     }}
  {\endlist}
  {\item}

\DeclareBibliographyDriver{techreport}{%
  \printnames{author}%
  \setunit{\labelnamepunct}\newblock
  \printfield[title]{title}%
  \newunit
  \printfield{institution}%
  \iffieldundef{number}{}{\setunit{\addcomma\space}\printfield{number}}%
  \iffieldundef{pages}{}{\setunit{\addcomma\space}\printfield{pages}}%
  \setunit{\addcomma\space}\printfield{year}%
  \newunit\newblock
}

\usepackage[english]{babel}
\usepackage{stackengine}

\definecolor{dkred}{rgb}{0.8,0,0}
\definecolor{blue}{rgb}{0,0,1}

\def\BibTeX{{\rm B\kern-.05em{\sc i\kern-.025em b}\kern-.08em
    T\kern-.1667em\lower.7ex\hbox{E}\kern-.125emX}}
    
\begin{document}

\title{Voltage-Aware Grid Aggregation: Expanding the European High-Voltage Network \\

\thanks{Funded by the European Union (ERC, NetZero-Opt, 101116212). Views and opinions expressed are however those of the author(s) only and do not necessarily reflect those of the European Union or the European Research Council. Neither the European Union nor the granting authority can be held responsible for them.}
}

\author{\IEEEauthorblockN{Benjamin Stöckl, Marco Anarmo, Sonja Wogrin and Yannick Werner}
\IEEEauthorblockA{\textit{Institute of Electricity Economics and Energy Innovation, Graz University of Technology }\\
\textit{Research Center ENERGETIC, Graz University of Technology}\\
Graz, Austria\\
\{benjamin.stoeckl,marnaizmontero,wogrin,yannick.werner\}@tugraz.at}}

\maketitle

\newcommand\copyrighttext{%
  \footnotesize \textcopyright \the\year{} IEEE. Personal use of this material is permitted. Permission from IEEE must be obtained for all other uses, including reprinting/republishing this material for advertising or promotional purposes, collecting new collected works for resale or redistribution to servers or lists, or reuse of any copyrighted component of this work in other works.}

\newcommand\copyrightnotice{%
\begin{tikzpicture}[remember picture,overlay]
\node[anchor=south,yshift=10pt] at (current page.south) {\fbox{\parbox{\dimexpr0.75\textwidth-\fboxsep-\fboxrule\relax}{\copyrighttext}}};
\end{tikzpicture}%
}

\newcommand\IEEEcopyrighttext{\textbf{979-8-3195-3554-2/26/\$31.00 ©2026 IEEE}}
\newcommand\IEEEcopyrightnotice{%
\begin{tikzpicture}[remember picture,overlay]
\node[anchor=south west, xshift=1.6cm, yshift=1.3cm] 
at (current page.south west) {%
 \IEEEcopyrighttext};
\end{tikzpicture}%
}

\copyrightnotice 
\vspace{-10pt}

\begin{abstract}
Energy system optimization models are indispensable for planning the European energy transition. Yet their applicability is constrained by the fundamental trade-off between spatial detail and computational tractability. Modelers often tackle this by spatially aggregating electricity networks. Existing methods, however, neglect differences in voltage levels, reducing them to a single level and thereby overlooking the critical role of transformers in expansion planning. Therefore, we propose a novel voltage-aware network partitioning and aggregation methodology that preserves individual voltage levels and transformers. We demonstrate the effectiveness of this approach and compare it against a voltage-unaware grid aggregation by solving a network expansion problem for a European case study using PyPSA. Our findings show that the proposed methodology preserves up to \qty{70}{\percent} of the transformer expansion costs in the aggregated model compared to the full grid model, thereby significantly improving the accuracy of investment decisions for transformers in the aggregated grid.
\end{abstract}

\begin{IEEEkeywords}
Spatial aggregation, network partitioning, transformers, DC power flow, expansion planning, clustering
\end{IEEEkeywords}

\section{Introduction}
\label{s:intro}

Due to decreasing losses at higher voltages, multiple voltage levels have evolved for transmitting electricity over long distances in real-world power grids.
In Europe, for example, \qtylist{220; 380}{\kilo\volt} are most common, with around \num{135} and \num{175} thousand \unit{km} line length, respectively~\cite{noauthor_statistical_2025}.
The lower losses at higher voltage levels, however, come at the cost of significantly increased line investment costs for better insulators and increased material usage~\cite{noauthor_comparison_2025}. Additionally, transferring electricity between voltage levels requires transformers, which entail substantial investment costs~\cite {noauthor_netzentwicklungsplan_2021}.
Consequently, the trade-off between network expansion on different voltage levels in real-world electricity grids is generally guided by the geographical transmission distance and associated losses~\cite{oeding_elektrische_nodate}. 

In recent years, energy system optimization models (ESOMs) have become a popular tool for analyzing this trade-off and supporting decision-making in practice. The complexity of real-world electricity grids, however, often renders those large-scale ESOMs computationally intractable. To overcome this, modelers simplify the physics of power flows or employ aggregation methods to reduce the spatial dimensions \cite{patil_advanced_2022, hoffmann_review_2024}. Some ESOMs, such as AnyMOD~\cite{goke_graph-based_2021}, EMPIRE~\cite{backe_empire_2022}, ETHOS.FINE~\cite{klutz_ethosfine_2025}, and Switch~\cite{johnston_switch_2019} model power flow as a transportation problem, while others, like GenX~\cite{bonaldo_2024_GenX} and LEGO~\cite{Wogrin2022}, implement DC power flow but lack spatial aggregation capabilities. To the best of the authors' knowledge, only PyPSA~\cite{brown_pypsa_2018}, Tulipa~\cite{siqueira_tulipa_nodate} and Spine~\cite{kouveliotis-lysikatos_network_2020} combine DC power flow with spatial aggregation. We refer the interested reader to~\cite{javanmardi_unraveling_2025} for a comprehensive review of how models manage spatial complexity.

Over the last few years, several methods to improve the partitioning process for DC optimal power flow applications have been investigated in the literature, including renewable energy resource potentials and electrical distance measures~\cite{phillips_spatial_2023, akdemir_open-source_2024}. Also, the use of locational marginal prices was examined in many studies to identify aggregations that account for line bottlenecks or form bidding zones with similar marginal prices~\cite{cao_incorporating_2018, chicco_overview_2019}. In a recent paper, we introduced a novel partitioning metric that is sensitive to line congestions~\cite{stockl_congestion-sensitive_2025}.

Although extensive literature explores distance metrics for partitioning, these approaches typically ignore voltage levels. Consequently, researchers often accept the loss of structural information regarding voltage levels and transformer infrastructure during aggregation \cite{frysztacki_comparison_2022, cao_incorporating_2018, phillips_spatial_2023, pache_e-highway_2015}. 

So far, no ESOM used for large-scale energy system analysis has employed a spatial aggregation technique that preserves voltage levels and thereby the critical role of transformers in real power systems. This often impedes the accurate assignment of individual line expansion costs across voltage levels and leads to capturing only a tiny fraction of the transformers in the original grid. As a consequence, ESOM outcomes tend to underestimate the required investment to reinforce transformers and provide little guidance on which voltage-level line expansion should occur. Both aspects are crucial for the expansion of real-world high-voltage electricity networks.

In this paper, we overcome these gaps by proposing a novel \emph{voltage-aware} network partitioning strategy, which preserves 1) individual voltage levels present in the original grid and 2) transformers as critical electricity infrastructure elements. We then apply this strategy to a case study of the European power system and demonstrate, using PyPSA, how it changes the investment requirements for lines and transformers compared to state-of-the-art \emph{voltage-unaware} approaches. To this end, we have implemented all methods in the standalone Network Partitioning and Aggregation Package (NPAP)~\cite{anarmo_npap_2026} in Python, ensuring reproducibility and extendability to any possible ESOM framework.

The remainder of this paper is organized as follows. Section~\ref{s:voltage-aware_grid_aggregation} introduces the proposed voltage-aware network partitioning and aggregation methodology and highlights its differences to conventional approaches using an illustrative example. In Section~\ref{s:european_case_study}, we provide a detailed description of the European case study and how we analyze line and transformer investment requirements using PyPSA. Section~\ref{s:results} presents and discusses the results of the grid aggregation process and investment outcomes, before Section~\ref{s:conclusion} concludes the paper.
%
\section{Towards voltage-aware grid aggregation}
\label{s:voltage-aware_grid_aggregation}

After introducing the network partitioning and aggregation methods in~\ref{s:network_partitioning_and_aggregation}, we show a stylized example in Section~\ref{s:illustrative_example}.


\subsection{Network partitioning and aggregation} 
\label{s:network_partitioning_and_aggregation}
Let $i,j,n \in \mathcal{N}$ denote the set of nodes and $v \in \mathcal{V}$ the set of voltage levels. Then the set of nodes (buses) can be divided based on the individual voltage level of each node~$n$, such that each subset $\mathcal{N}_v \subseteq \mathcal{N}$ contains all nodes~$n$ that belong to voltage level~$v$. Then we define \emph{network partitioning} as the process of finding a mapping $\mathcal{M}$ from the set of nodes $\mathcal{N}$ to the set $k \in \mathcal{K}$ of disjoint subsets or clusters. We use $\mathcal{N}_k \subseteq \mathcal{N}$ with $\cup_{k \in \mathcal{K}} \mathcal{N}_k = \mathcal{N}$, to denote the subset of nodes belonging to cluster~$k$, such that each node $n$ belongs to exactly one cluster~$k$.

The common approach in ESOMs has so far been to find a partition of the network regardless of the nodes' voltage levels~\cite{hoffmann_review_2024}. We refer to this as \emph{voltage-unaware (VU) partitioning}. Formally, the goal is to find a mapping $\mathcal{M}^{\mathrm{VU}}$, which maps the set of nodes $\mathcal{N}$ to the set of clusters $\mathcal{K}$, i.e., $\mathcal{M}^{\mathrm{VU}}: \mathcal{N} \mapsto \mathcal{K}$. In contrast, we propose a novel \emph{voltage-aware (VA) partitioning} that explicitly accounts for the voltage level~$v_n$ of each node~$n$. To achieve that, we require that all nodes within a single cluster must have the same voltage level. Let $\mathcal{K}_v$ denote the subset of clusters $\mathcal{K}$ that only contain nodes on voltage level~$v$. Then we define mapping $\mathcal{M}^{\mathrm{VA}} : (\mathcal{N}_v \mapsto \mathcal{K}_v)_{\forall v \in \mathcal{V}}$. We provide additional information on the practical implementation of the VA partitioning in Appendix~\ref{a:distances}.

We use the term \emph{network aggregation} to describe the process of aggregating a network according to a given partitioning (mapping). In this paper, we follow a copperplate approach that, for each cluster $k$, combines all buses within this cluster, $n \in \mathcal{N}_k$, into a single node. All generators, loads, etc., connected to one of the buses are aggregated accordingly, such that, e.g., installed capacities are summed, and renewable capacity factors are averaged~\cite{frysztacki_comparison_2022}. All lines connecting nodes within a cluster are removed. Additionally, we aggregate parallel lines between distinct clusters $k$ and $k'$ by summing their capacity and deriving an equivalent reactance using standard calculus. Importantly, we also derive a capacity-weighted average of line lengths and investment costs for use in ESOM expansion planning later. Because transformers are effectively represented as lines in DC power flow, the same aggregation strategy applies to them.
In the following, we use the term \emph{grid aggregation} to refer to the combined process of network partitioning and aggregation.
\subsection{Illustrative example}
\label{s:illustrative_example}
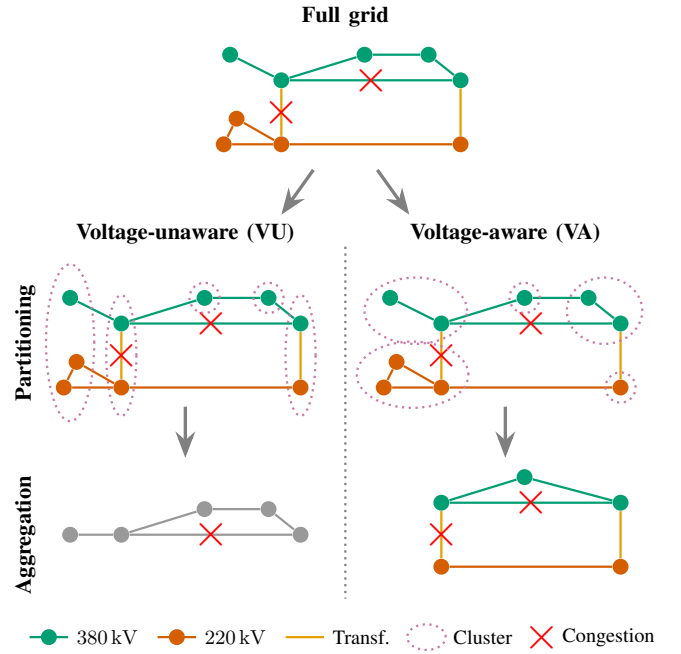
\begin{figure}[bt]
    \centering
    \resizebox{0.5\textwidth}{!}{
        \begin{tikzpicture}[
    node distance=2cm,
    bus/.style={circle, fill, inner sep=2.5pt},
    bus380/.style={bus, color=oiBluishGreen},
    bus220/.style={bus, color=oiVermilion},
    busAgg/.style={bus, color=gray!80},
    line380/.style={oiBluishGreen, line width=1pt}, 
    line220/.style={oiVermilion, line width=1pt}, 
    lineTrafo/.style={oiOrange, line width=1pt},
    lineAgg/.style={gray!80, line width=1pt},
    cluster/.style={ellipse, draw=oiReddishPurple, dotted, line width=1pt, inner sep=1pt},
    arrow/.style={-{Stealth[scale=1.2]}, thick, gray},
    congested/.style={midway, text=red, font=\huge,}
]

\newcommand{\drawgrid}[8]{
    \node (#3R1) at (#1+0.2, #2+1.2) [#4, color=oiBluishGreen] {};
    \node (#3R2) at (#1+1  , #2+0.8) [#4, color=oiBluishGreen] {};
    \node (#3R3) at (#1+2.3, #2+1.2) [#4, color=oiBluishGreen] {};
    \node (#3R4) at (#1+3.3, #2+1.2) [#4, color=oiBluishGreen] {};
    \node (#3R5) at (#1+3.8, #2+0.8) [#4, color=oiBluishGreen] {};
    
    \draw [#5] (#3R1) -- (#3R2);
    \draw [#5] (#3R2) -- (#3R3);
    \draw [#5] (#3R3) -- (#3R4);
    \draw [#5] (#3R4) -- (#3R5);
    \draw [#5] (#3R2) -- node[#8] {$\times$} (#3R5);
    
    \node (#3G1) at (#1+0.3, #2+0.2) [#4, color=oiVermilion] {};
    \node (#3G2) at (#1+0.1, #2-0.2) [#4, color=oiVermilion] {};
    \node (#3G3) at (#1+1  , #2-0.2) [#4, color=oiVermilion] {};
    \node (#3G4) at (#1+3.8, #2-0.2) [#4, color=oiVermilion] {};
    
    \draw [#6] (#3G1) -- (#3G2);
    \draw [#6] (#3G2) -- (#3G3);
    \draw [#6] (#3G3) -- (#3G4);
    \draw [#6] (#3G1) -- (#3G3);
    
    \draw [#7] (#3R2) --node[#8] {$\times$} (#3G3);
    \draw [#7] (#3R5) -- (#3G4);
}

\begin{scope}[yshift=5.3cm, xshift=0cm]
    \node[font=\bfseries] at (2, 1.8) {Full grid};
    \drawgrid{0}{0}{init}{bus}{line380}{line220}{lineTrafo}{congested}
\end{scope}

\draw[arrow, line width = 1.5pt] (1.5, 4.7) -- (1, 4);
\draw[arrow, line width = 1.5pt] (2.5, 4.7) -- (3, 4);

\node[font=\bfseries] at (-0.5, 3.7) {Voltage-unaware (VU)};
\node[font=\bfseries] at (4.5, 3.7) {Voltage-aware (VA)};

\node[font=\bfseries, rotate=90, align=center] at (-3, 2) {Partitioning};
\node[font=\bfseries, rotate=90, align=center] at (-3, -1) {Aggregation};

\begin{scope}[yshift=1.5cm, xshift=0cm]
    \drawgrid{-2.5}{0}{vuP}{bus}{line380}{line220}{lineTrafo}{congested}
    \node[cluster, fit=(vuPR1) (vuPG1) (vuPG2)] {};
    \node[cluster, fit=(vuPR2) (vuPG3)] {};
    \node[cluster, fit=(vuPR3)] {};
    \node[cluster, fit=(vuPR4)] {};
    \node[cluster, fit=(vuPG4) (vuPR5)] {};
\end{scope}

\begin{scope}[yshift=1.5cm, xshift=2.5cm]
    \drawgrid{0}{0}{vaP}{bus}{line380}{line220}{lineTrafo}{congested}
    \node[cluster, fit=(vaPR1) (vaPR2)] {};
    \node[cluster, fit=(vaPR4) (vaPR5)] {};
    \node[cluster, fit=(vaPG1) (vaPG2) (vaPG3)] {};
    \node[cluster, fit=(vaPG4)] {};
    \node[cluster, fit=(vaPR3)] {};
\end{scope}

\draw[arrow, line width = 1.5pt] (-0.5, 1) -- (-0.5, 0.2);
\draw[arrow, line width = 1.5pt] (4.5, 1) -- (4.5, 0.2);

 \draw[gray, dotted, line width=1pt] (2, 3.5) -- (2, -2);

\begin{scope}[yshift=-1.5cm, xshift=-2.5cm]
    \node (vA1) at (0.2  , 0.5) [busAgg] {};
    \node (vA2) at (1  , 0.5) [busAgg] {};
    \node (vA3) at (2.3, 0.9) [busAgg] {};
    \node (vA4) at (3.3, 0.9) [busAgg] {};
    \node (vA5) at (3.8, 0.5) [busAgg] {};
    \draw[lineAgg] (vA1) -- (vA2);
    \draw[lineAgg] (vA2) -- (vA3);
    \draw[lineAgg] (vA2) --node[midway, text=red, font=\huge, ] {$\times$} (vA5);
    \draw[lineAgg] (vA3) -- (vA4);
    \draw[lineAgg] (vA4) -- (vA5);
\end{scope}

\begin{scope}[yshift=-1.5cm, xshift=2.5cm]
    \node (vaAR1) at (1, 1) [bus380] {};
    \node (vaAR2) at (3.8, 1) [bus380] {};
    \node (vaAR3) at (2.3, 1.4) [bus380] {};
    \draw[line380] (vaAR1) --node[congested ] {$\times$} (vaAR2);
    \draw[line380] (vaAR1) -- (vaAR3);
    \draw[line380] (vaAR3) -- (vaAR2);
    
    \node (vaAG1) at (1, 0) [bus220] {};
    \node (vaAG2) at (3.8, 0) [bus220] {};
    \draw[line220] (vaAG1) -- (vaAG2);

    \draw[lineTrafo] (vaAR1) --node[congested] {$\times$} (vaAG1);
    \draw[lineTrafo] (vaAR2) -- (vaAG2);
\end{scope}

\node[anchor=center] at (2, -2.6) {
    \begin{tikzpicture}
        \draw[line380] (0, 0) --node[midway, bus380] {} (0.6, 0) node[right, font=\small, text=black] {\qty{380}{\kilo\volt}};
        
        \draw[line220] (2, 0) -- node[midway, bus220] {} (2.6, 0) node[right, font=\small, text=black] {\qty{220}{\kilo\volt}};
    
        \draw[lineTrafo] (4, 0) -- (4.6, 0) node[right, font=\small, text=black] {Transf.};
        
        \node[cluster, minimum width=0.6cm, minimum height=0.4cm] (legCluster) at (6.2, 0) {};
        \node[right, font=\small, text=black] at (6.5, 0) {Cluster};
        
        \node[text=red, font=\huge] at (8, 0)  {$\times$};
        \node[right, font=\small, text=black] at (8.2, 0) {Congestion};
    \end{tikzpicture}
    };

\end{tikzpicture}
    }
    \caption{Illustration of voltage-unaware and voltage-aware network partitioning and aggregation.}
    \label{fig:stylized_testcase}
    \vspace{-0.6cm}
\end{figure}

We illustrate the functionality of voltage-aware partitioning and highlight its differences to voltage-unaware partitioning through an illustrative example depicted in Fig.~\ref{fig:stylized_testcase}.
The top part of the figure shows a stylized high-voltage grid with two exemplary voltage levels of \qty{220}{\kilo\volt} (orange) and \qty{380}{\kilo\volt} (green) lines, connected by two transformers (yellow), shown in a side-view for visual clarity.
We now assume that one of the \qty{380}{\kilo\volt} lines and one transformer are frequently congested (red cross), and that expanding their capacity may be worth considering.

In the row below, we illustrate a potential state-of-the-art voltage-unaware (left) and proposed voltage-aware (right) partitioning. The cluster assignment of nodes is indicated by dashed circles.
Standard VU partitioning methods tend to assign nodes to a single cluster, regardless of their distinct voltage levels. When aggregating the network using this partitioning, the process fails to preserve the difference in voltage levels and frequently removes transformers. Hence, the VU grid aggregation retains only the congested line, without preserving voltage levels, potentially failing to capture accurate line expansion cost and the investment required to reinforce the transformers.

In contrast, when using the VA partitioning, shown on the right, clusters are formed only within their respective voltage levels. During the aggregation, this preserves the two separate voltage levels and the transformers between them. In particular, it allows the capture of congestion and investment needs for transformers, which are now, by construction, present in the aggregated grid. Additionally, the representation of individual voltage levels enables the use of line investment costs for lines at different voltage levels. Both aspects foster a more accurate representation of the physical infrastructure and associated expansion costs in the ESOMs, improving decision support.

\section{Expanding Europe's high-voltage network}
\label{s:european_case_study}
After introducing the novel concept of VA partitioning, we now want to demonstrate the idea using a case study of the European power system and the ESOM PyPSA~\cite{brown_pypsa_2018}. In Section~\ref{s:case_study}, we provide a detailed description of the network topology used in the case study. Afterwards, we briefly elaborate in Section~\ref{s:PyPSA_description} how we integrated everything in PyPSA and provide an overview of the configuration settings. 

\subsection{Voltage-aware grid representation}
\label{s:case_study}
\begin{figure}
    \centering
    \includegraphics[width=1\linewidth]{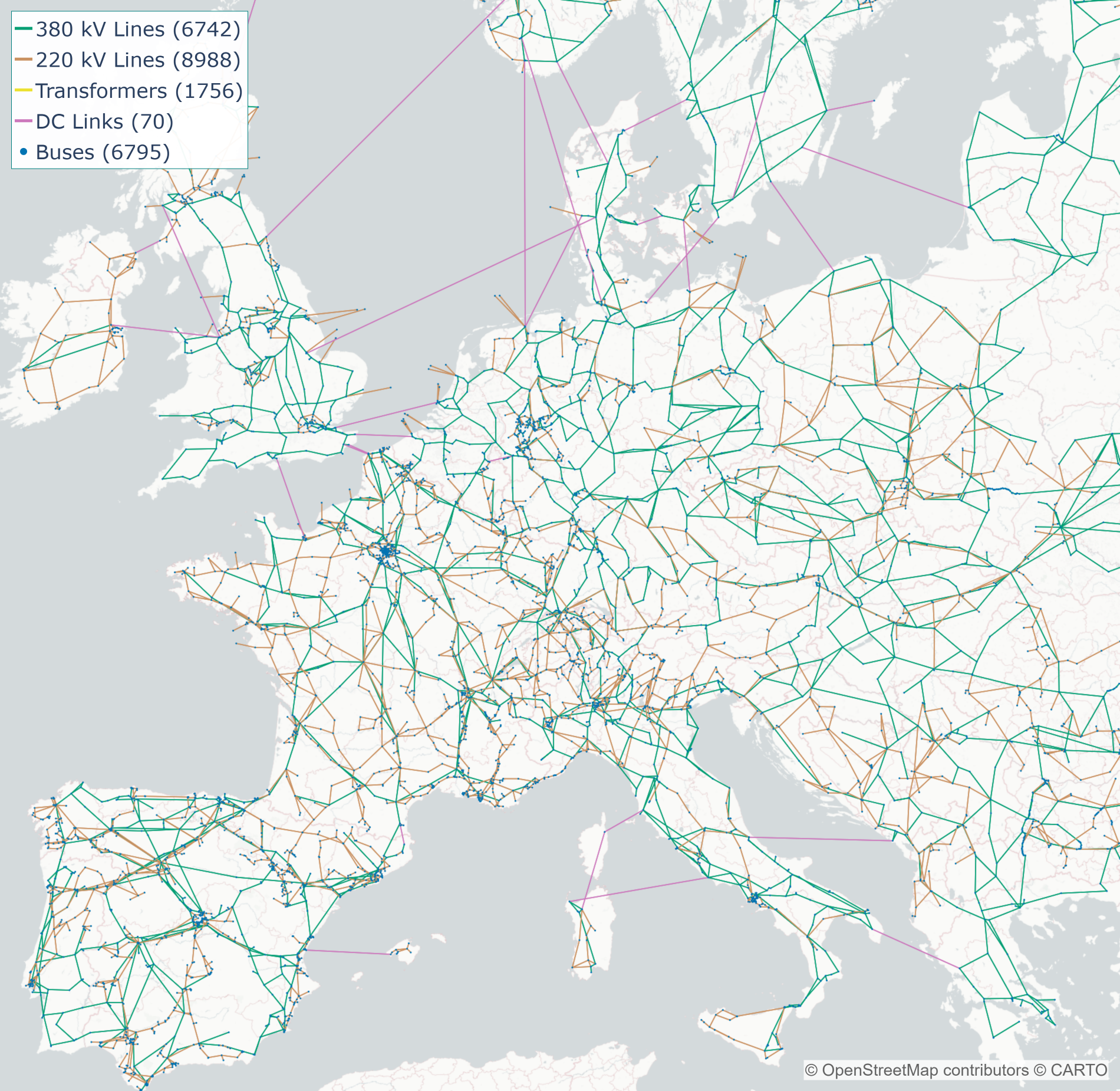}
    \caption{Representation of the European high-voltage network taken from \cite{Xiong2025} version 0.7.0 \cite{Xiong2026} using NPAP \cite{anarmo_network_nodate}.}
    \label{fig:Europe_full_grid}
    \vspace{-0.6cm}
\end{figure}

In Fig.~\ref{fig:Europe_full_grid}, we show the topology of the European high-voltage network taken from~\cite{Xiong2025} in its most recent version~\cite{Xiong2026} and illustrated using NPAP \cite{anarmo_network_nodate}. 
We map all power lines to \qty{220}{\kilo\volt} and \qty{380}{\kilo\volt} voltage levels and do not aggregate stub lines. The network of the full model contains \num{6795}~buses, connected by \num{8988} \qty{220}{\kilo\volt} lines with a total length of around \num{110} thousand \unit{\kilo\metre} and \num{6742} \qty{380}{\kilo\volt} lines with a total length of around \num{160} thousand \unit{\kilo\metre}, as well as \qty{1756}{} transformers.

A brief overview of the reinforcement investment costs for lines and transformers at their respective voltage levels is shown in Table~\ref{tab:expansion_costs}.
To account for the different roles and technical characteristics of \qty{220}{\kilo\volt} and \qty{380}{\kilo\volt} lines discussed in Section~\ref{s:intro}, we differentiate between their investment cost. Analyzing the per-kilometer costs in Table~\ref{tab:expansion_costs}, \qty{380}{\kilo\volt} lines are, on average, 1.6 times more expensive to build than \qty{220}{\kilo\volt} lines. However, because their transmission capacity is approximately three times greater, the specific investment cost per unit of capacity (\si{\mega\volt\ampere/\kilo\metre}) is lower for the higher voltage level. Despite this economy of scale, using \qty{380}{\kilo\volt} infrastructure requires power transformers, incurring substantial additional investment costs.
As detailed in Table~\ref{tab:expansion_costs}, the capital required to bridge distinct voltage layers is substantial, with a single \qty{380}{}/\qty{220}{\kilo\volt} power transformer costing between \qty{4.6}{} and \qty{8.5}{\mega\euro}, not including associated facility assets, such as switchgear. 
For our case study, we use the technology-specific average investment cost and note again that we account only for reinforcement, rather than new construction, which leads to comparably low costs in the ESOM.

While the data in~\cite{Xiong2025} provide good estimates of line lengths and capacities, we found that the data on transformer capacities and reactances can be improved. To do so, we geographically match transformers in~\cite{Xiong2025} to those in the JAO Core Static Grid Model (CSGM)~\cite{joint_allocatoin_office_static_nodate} using a two-step procedure. First, we extract substation names from the CSGM and geocode them using the Nominatim API~\cite{nominatim_nominatim_nodate} and the Overpass API~\cite{overpass_api_overpass_nodate}. Second, we match entries based on their geographical distance, their connected voltage levels, and the degree of name matching. 
Of the 495 PyPSA transformers (i.e., substations) located in countries covered by the CSGM data set, 195 were successfully matched (\qty{39}{\percent}), corresponding to 408 CSGM entries, out of 520 in total. For the matched entries, we compare the capacities in the original PyPSA dataset~\cite{Xiong2025} and those in CSGM~\cite{joint_allocatoin_office_static_nodate}, and find that, on average, they are roughly overestimated by a factor of approximately~4. 
For the 683 out of 878 (\qty{78}{\percent}) substations in the PyPSA data set that were unmatched, we replace the capacity and reactance with the median of the voltage-pair values for the matched CSGM substations.

\begin{table}
    \caption{Reinforcement investment costs for \qtylist{220; 380; 400}{\kilo\volt} lines with two circuits and transformers connecting \setlength{\tabcolsep}{4pt}\qtylist{380; 220}{\kilo\volt}.}
    \label{tab:expansion_costs}
    \centering
    \begin{tabular}{lcccrc}
        \toprule
        \textbf{Type} & \shortstack{\textbf{Voltage} \\ in~\si{\kilo\volt}} & \shortstack{\textbf{Cost} \\ in \si{\mega\euro\per(\kilo\metre)}} & \shortstack{\textbf{Cost} \\ in \si{\mega\euro}} & \shortstack{\textbf{Capacity} \\ in \si{\mega\volt\ampere}} & \textbf{Source} \\
        \midrule
         Line    & 220    &  0.44 -- 0.53  & --    &  983.1     &   \cite{noauthor_unit_2025}     \\
         Line    & 380    & 1.26   & --   &  3396.2    &    \cite{noauthor_unit_2025}    \\
         Line    & 380    & 0.70   & --    &  1698.1    &   \cite{noauthor_netzentwicklungsplan_2021}     \\
        Line    & 400    & 0.67    & --   &  3574.9    &  \cite{noauthor_unit_2023}      \\
        \midrule
         Transf.  & 380/220 & --    &  4.6       & 600.0          &  \cite{noauthor_unit_2025}      \\
         Transf.    & 380/220 & --    &  8.5       & 600.0           &  \cite{noauthor_netzentwicklungsplan_2021} \\
        \bottomrule
    \end{tabular}
    \vspace{-0.4cm}
\end{table}
\subsection{Voltage-aware expansion planning with PyPSA}
\label{s:PyPSA_description}
In the following, we briefly explain how we prepared the data for PyPSA~\cite{brown_pypsa_2018}, the configuration used, and how we implemented our VA grid aggregation workflow. The data and full implementation of the workflow are available on GitHub~\cite{stockl_voltage-aware_2026}.
We retrieve the data for generators and loads using the PyPSA-Eur model-building pipeline, but we skip the step of spatial aggregation~\cite{horsch_pypsa-eur_2023} to preserve voltage levels and only consider the power system. We configure the case study with a 2050 planning horizon and use 2013 as the historical weather year. We further temporally aggregate the 8760 hourly time steps into three representative 24-hour periods using k-means clustering~\cite{tsam_2022}.

Before running the transmission expansion planning problem, we first solve an operational problem across the entire grid to identify investment candidates for lines and transformers. We classify each component as an investment candidate if its loading exceeds either \qty{70}{\percent} across all time steps or \qty{90}{\percent} at least once. Subsequently, we perform the VU and VA network partitioning and aggregation, described in Section~\ref{s:network_partitioning_and_aggregation}, based on k-medoids clustering and geographical distances using NPAP's \texttt{geographical\_kmedoids\_haversine} partitioning strategy and \texttt{GEOGRAPHICAL} aggregation mode applicable for PyPSA. For that, we have integrated NPAP directly into the PyPSA-Eur workflow~\cite{anarmo_integrate_nodate}. Finally, we solve the continuous transmission expansion planning problem for the full, VU, and VA aggregated networks, using the Barrier method without crossover in Gurobi~v. 13.0.0 \cite{Gurobi2023}.
\section{Results}
\label{s:results}
In the following, we first compare in Section~\ref{s:results_part} the VU and VA partitioning and aggregation results for the European case study introduced in~\ref{s:european_case_study}. Afterwards, in Section~\ref{s:res_exp_plan}, we analyze, using PyPSA, how the different topologies lead to different optimal line and transformer expansion decisions by comparing them with the model outcomes for the full grid.
\subsection{Partitioning and aggregation}
\label{s:results_part}
\begin{figure}
    \centering
    \resizebox{\columnwidth}{!}{
        \input{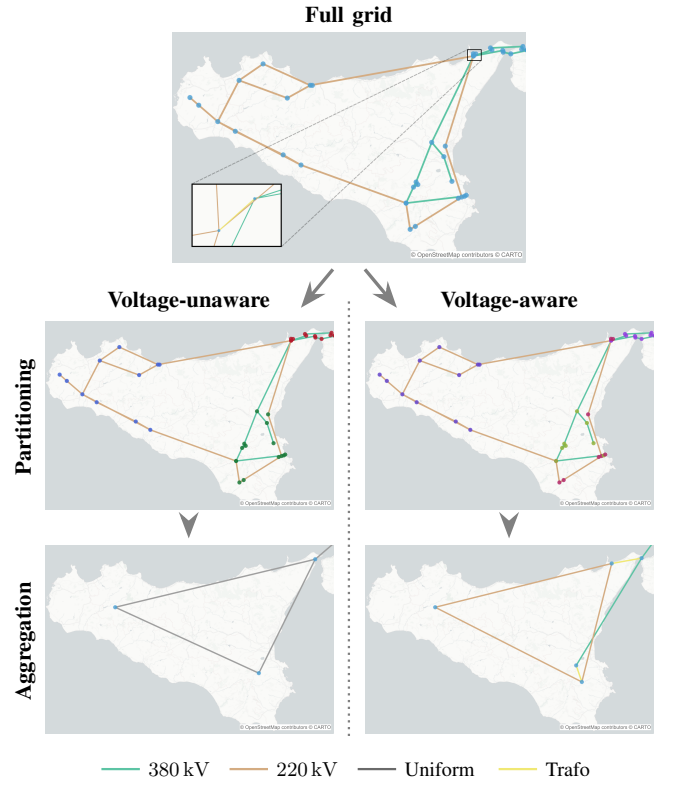}
    }
    \caption{Network partitioning and aggregation results obtained with NPAP \cite{anarmo_network_nodate} for a close-up of Sicily. Circle colors in the partitioning indicate cluster assignment.}
    \label{fig:results_npap}
    \vspace{-0.5cm}
\end{figure}

In Fig.~\ref{fig:results_npap}, we zoom into Sicily and show the results for the state-of-the-art VU and the novel VA network partitioning and aggregation obtained with NPAP~\cite{anarmo_network_nodate} with a total of $k=500$ nodes in the aggregated grids.
The subfigure in the upper part shows the original topology, with a substation in the northern part of the island magnified, hosting a transformer connecting the \qtylist{220;380}{\kilo\volt} voltage levels. Below, on the left, the VU aggregation yields a grid topology with three buses distributed evenly geographically and without transformers. In contrast, the VA grid aggregation preserves individual voltage levels and transformers by grouping nodes into five clusters. While this preserves essential topological features, it highlights a drawback of the proposed VA approach: it tends to place buses at different voltage levels that are geographically close, resulting in sparser coverage elsewhere. Interestingly, across Europe, the VU grid aggregation with $k=500$ nodes contains \num{955} connecting lines but no transformers. This occurs due to the geographical proximity of low- and high-voltage busbars within substations. For the same number of clusters $k$, the VA grid aggregation features fewer lines (\num{754}) but preserves \num{344} transformers.
%
\subsection{Expansion planning}
\label{s:res_exp_plan}
After analyzing the different topologies between VU and VA grid aggregation, we now examine the optimal expansion planning decisions obtained with PyPSA~\cite{brown_pypsa_2018}. In the following, we refer to the optimization model representing the full grid (FG) as the FG model and to the models representing the VU and VA aggregations as the VU and VA models, respectively.

\begin{figure*}
    \centering

    \begin{tikzpicture}[
        node distance=2cm,
        bus/.style={circle, fill, inner sep=2.5pt},
        bus380/.style={bus, color=oiVermilion},
        bus220/.style={bus, color=oiBluishGreen},
        busAgg/.style={bus, color=gray!80},
        line380/.style={oiVermilion, line width=1pt},
        line220/.style={oiBluishGreen, line width=1pt},
        lineTrafo/.style={oiOrange, line width=1pt},
        lineAgg/.style={gray!80, line width=1pt},
        cluster/.style={ellipse, draw=oiReddishPurple, dotted, line width=1pt, inner sep=1pt},
        arrow/.style={-{Stealth[scale=1.2]}, thick, gray},
        congested/.style={midway, text=red, font=\huge,}
    ]

        \draw[NPAP_380, line width=1pt] (0, 0) -- (0.6, 0) node[right, font=\small, text=black] {\qty{380}{\kilo\volt}}; 
        \draw[NPAP_380_inv, line width=1.5pt] (2, 0) -- (2.6, 0) node[right, font=\small, text=black] {\qty{380}{\kilo\volt} Inv.}; 
        
        \draw[NPAP_220, line width=1pt] (4.6, 0) -- (5.2, 0) node[right, font=\small, text=black] {\qty{220}{\kilo\volt}};   
        \draw[NPAP_220_inv, line width=1.5pt] (6.6, 0) -- (7.2, 0) node[right, font=\small, text=black] {\qty{220}{\kilo\volt} Inv.};   
        
        \draw[NPAP_uniform, line width=1pt] (9.2, 0) -- (9.8, 0) node[right, font=\small, text=black] {Uniform};
        
        \draw[NPAP_transformer, line width=1pt] (11.4, 0) -- (12, 0) node[right, font=\small, text=black] {Trafo};
        \draw[NPAP_transformer_inv, line width=1.5pt] (13.2, 0) -- (13.8, 0) node[right, font=\small, text=black] {Trafo Inv.};
        
    \end{tikzpicture}\\ 
    \vspace{-0.2cm}

        \subfloat[Full grid.\label{fig:map_expanded_lines_full_grid}]{%
            \includegraphics[width=0.32\linewidth]{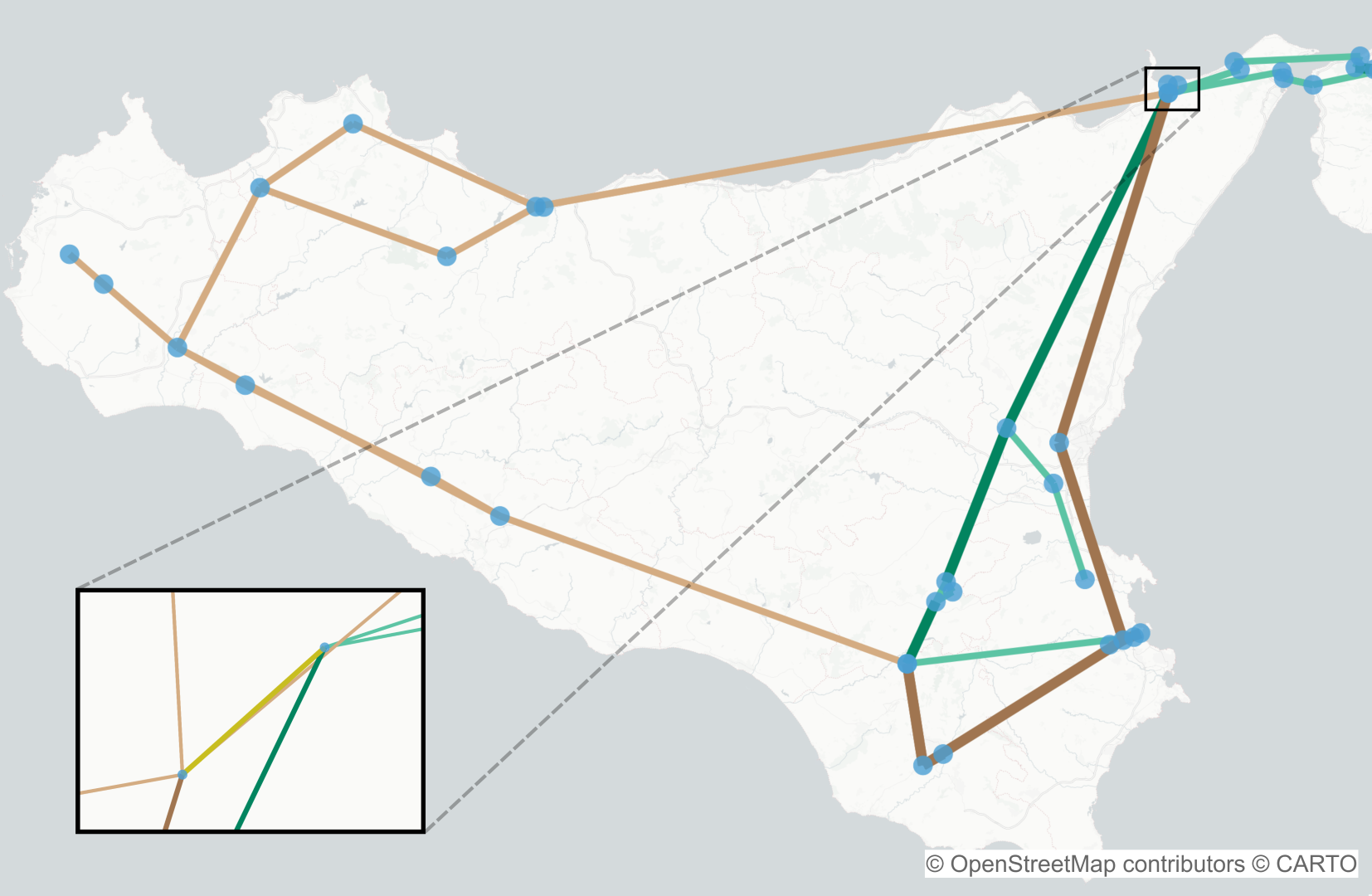}%
        }
        \hfil 
        \subfloat[Voltage-unaware aggregation.\label{fig:map_expanded_lines_VU}]{%
            \includegraphics[width=0.32\linewidth]{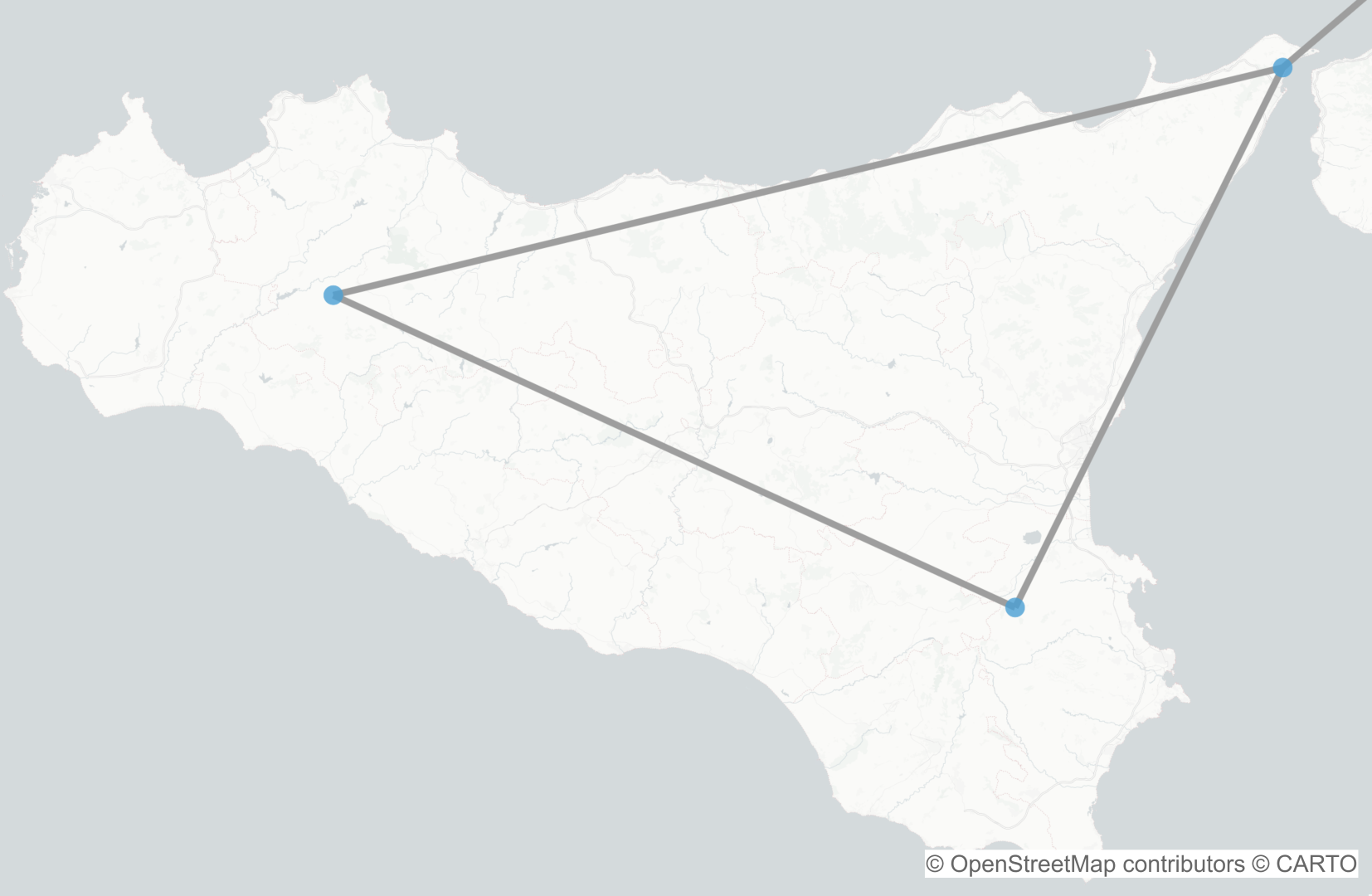}%
        }
        \hfil
        \subfloat[Voltage-aware aggregation.\label{fig:map_expanded_lines_VA}]{%
            \includegraphics[width=0.32\linewidth]{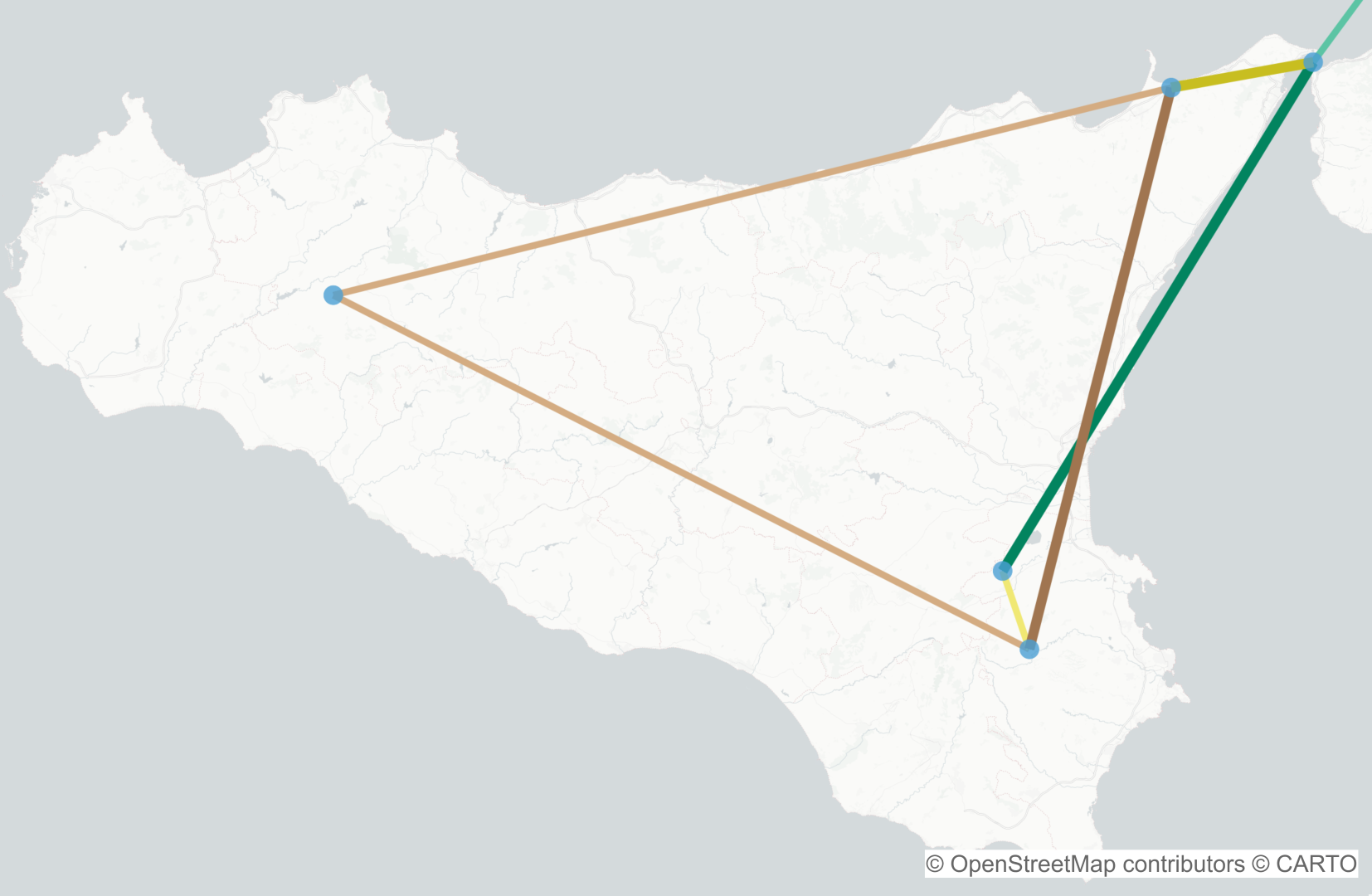}%
        }
    
    \caption{Network expansion results shown for Sicily for the models representing the full grid and the voltage-unaware and voltage-aware aggregations, respectively, with $k=500$~nodes. Darker line color and thickness indicate the capacity invested.}
  \label{fig:map_expanded_lines}
  \vspace{-0.4cm}
\end{figure*}
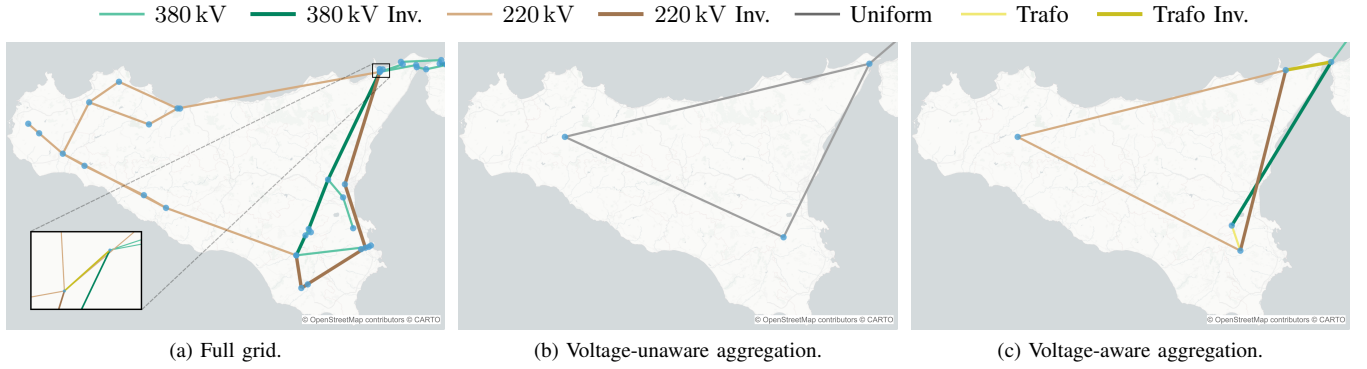

In Fig.~\ref{fig:map_expanded_lines}, we show a close-up of Sicily again and illustrate the expanded network capacities for the full FG, VU, and VA models with $k=500$~nodes. Darker line colors and thicker line widths indicate the amount of capacity expanded. In the FG model (Fig.~\ref{fig:map_expanded_lines_full_grid}), some \qty{220}{\kilo\volt} (orange) and \qty{380}{\kilo\volt} (green) lines, as well as one transformer (yellow) in the East, are expanded. This is not captured by the VU model, although the Eastern line is defined as an investment candidate (Fig.~\ref{fig:map_expanded_lines_VU}). The VA model, however, successfully captures the bottlenecks and expands both the \qtylist{220;380}{kV} lines (Fig.~\ref{fig:map_expanded_lines_VA}) and the transformer.
Next, we evaluate investment decisions for various numbers of nodes~$k$ in the aggregated models and compare them with the results from the FG model. Table~\ref{tab:investment_results} summarizes the total expanded capacity(-length), total investment costs, and cost deviations relative to the FG benchmark. Note that the shown costs are annualized and merely reflect the expansion of existing lines and transformers, and therefore do not reflect the actual cost of grid expansion in Europe.
\begin{table}[]
    \centering
    \caption{Optimal investment decisions, deviation costs, and work requirements for full grid (FG) and aggregated (VU, VA) models for various numbers of nodes~$k$.
    }
    \label{tab:investment_results}
    \setlength{\tabcolsep}{3pt}
    \begin{tabular}{ll rr rr r}
        \toprule
        \multirow{2}{*}{\textbf{Mod.}} & \multirow{2}{*}{\textbf{$k$}} & \multicolumn{2}{c}{\textbf{Lines}} & \multicolumn{2}{c}{\textbf{Transformers}} & \multirow{3}{*}{\textbf{\makecell{Work-\\units}}}\\
        \cmidrule(lr){3-4} \cmidrule(lr){5-6}
        & \text{(nodes)} & \textbf{Cap.-length} & \textbf{Cost (Dev.)} & \textbf{Capacity} & \textbf{Cost (Dev.)} & \\
        & & in \si{\giga\volt\ampere\kilo\meter} & in \si{\mega\euro} (\si{\percent}) & in \si{\giga\volt\ampere} & in \si{\mega\euro} (\si{\percent}) & \\
        \midrule
        \textbf{FG} & 6795 & \databarLine{3254.9} & 31.0 \phantom{0} (0.0) & \databarTrafo{13.4} & 1.8 (\phantom{0}0.0) & 1592.5\\
        \midrule
        \multirow{4}{*}{\textbf{VU}} 
        & 250  & \databarLine[oiVermilion]{1173.3} & 11.9 (-61.6) & \databarTrafo{0.0} & 0.0 (-100.0) &   35.0\\
        & 500  & \databarLine[oiVermilion]{1096.6} & 11.8 (-62.0) & \databarTrafo{0.0} & 0.0 (-100.0) &   97.4\\
        & 1000 & \databarLine[oiVermilion]{1212.0} & 13.4 (-56.9) & \databarTrafo{0.0} & 0.0 (-100.0) &  588.6\\
        & 3000 & \databarLine[oiVermilion]{1847.2} & 19.5 (-37.1) & \databarTrafo{0.0} & 0.0 (-100.0) & 1286.5\\
        \midrule
        \multirow{4}{*}{\textbf{VA}} 
        & 250  & \databarLine[oiBluishGreen]{343.6}  & 3.3  (-89.5) & \databarTrafo[oiBluishGreen]{0.6}  & 0.1 (-95.9) &  58.9 \\
        & 500  & \databarLine[oiBluishGreen]{1543.9} & 13.9 (-55.2) & \databarTrafo[oiBluishGreen]{9.4}  & 1.3 (-29.6) & 103.4 \\
        & 1000 & \databarLine[oiBluishGreen]{1096.8} & 10.6 (-66.0) & \databarTrafo[oiBluishGreen]{6.7}  & 0.9 (-51.1) & 415.1 \\
        & 3000 & \databarLine[oiBluishGreen]{2501.3} & 24.1 (-22.3) & \databarTrafo[oiBluishGreen]{8.8}  & 1.2 (-34.0) & 638.8 \\
        \bottomrule
    \end{tabular}
    \vspace{-0.5cm}
\end{table}
We first focus on line capacity expansion. Therefore, we compare the extended capacity multiplied by the line length, to compensate for the added capacity of a short route at a lower cost than a long route. Compared with the \qty{3254.9}{\giga\volt\ampere\kilo\meter} of total line expansion in the FG model, both aggregation methods underestimate the required capacity, even at higher node counts. The VA model with \qty{3000} nodes achieves the closest result, although it still underestimates the total line capacity expansion by \qty{23}{\percent}.
This suggests that grid aggregation generally leads to a significant loss of information about line investments, as many candidate lines are removed during aggregation.
For \num{500} and \num{1000} nodes, both VU and VA models show similar performance with capacities between \qtyrange{31.8}{33.7}{\giga\volt\ampere}.
Notably, for the coarsest aggregation ($k=250$), the VU approach underestimates the full model's capacity less than the VA approach does. We hypothesize that this is due to the VA aggregation splitting the number of clusters across two voltage levels, thereby yielding geographically large clusters, which leads to many candidates being removed. 
The line investment costs show a similar deviation than the epxanded capacity-length of typically underestimating the costs in the VU and VA models compared to the FG model.
 

We next compare the investment results for transformer expansion. In the FG model, transformers are expanded by around \qty{13.4}{\giga\volt\ampere} with an associated annualized cost of \qty{1.8}{\mega\euro}. In the VU models, there is no capacity expansion for transformers at all, as they are all removed during the aggregation process, as discussed in Section~\ref{s:results_part}. Consequently, using the VU grid aggregation provides little information or guidance on which transformers should be expanded.
In contrast, the VA models with \num{500} nodes or more are able to capture a large part of the transformer investment optimal in the FG model. This demonstrates that, by preserving individual voltage levels across the full grid, the proposed VA grid aggregation captures the essential role of transformers and their investment requirements in ESOMs. Surprisingly, both aggregations with \num{1000} and \num{3000} nodes yield slightly lower capacity expansion than the model with \num{500} nodes. This suggests that these specific configurations fail to preserve critical grid bottlenecks and highlights the drawback of nonhierarchical partitioning methods, in which cluster mappings can vary substantially as the number of clusters changes~\cite{stockl_congestion-sensitive_2025}.
\section{Conclusion}\label{s:conclusion}
A critical limitation of existing spatial aggregation techniques for ESOMs is the neglect of the different voltage levels present in real-world electricity grids, thereby overlooking differences in infrastructure expansion costs across voltage levels and the critical role of transformers. In this paper, we overcome this major drawback by proposing a novel VA grid aggregation methodology that preserves individual voltage levels during the network partitioning and aggregation process. We demonstrate and compare the proposed method with state-of-the-art VU grid aggregation using a realistic case study of the European high-voltage network, and show its effectiveness in preserving the topology of the original grid, particularly transformers. Using PyPSA~\cite{brown_pypsa_2018}, we solve a network expansion planning problem and show that the VA grid aggregation significantly improves the accuracy of transformer investment decisions by capturing critical bottlenecks that conventional techniques entirely obscure. Thereby providing better decision-making support for network expansion planning in real-world electricity grids.
Future research should explore how voltage-aware grid aggregation can be further improved by employing alternative network partitioning and aggregation methods, e.g., based on electrical distances or KRON-reduction. On the other hand, data on existing transformers should be further refined to obtain more realistic estimates of current installed capacities and potential investment requirements.


\appendices

\section{Distances}
\label{a:distances}
In the following, we briefly describe the implementation of the voltage-aware partitioning in NPAP. Let $\mathbf{D}\in \mathbb{R}^{n \times n}$ be the node distance matrix that captures, e.g., geographical proximity, and $p$ be the number of nodes. Then, to avoid clustering nodes with different voltage levels, we alter $\mathbf{D}$, such that for two voltage levels $v$ and $v'$, all nodes $n \in \mathcal{N}_v$ and $n' \in \mathcal{N}_{v'}$, get infinity distance assigned:
%
\begin{equation*}
\label{eq:distance_matrix}
\small
\setlength{\arraycolsep}{2.5pt}
 \mathbf{D} =
    \begin{pNiceArray}{ccc|ccc}[first-row, first-col, margin]
        & \Block{1-3}{\vspace*{8pt} n \in \mathcal{N}_{v}} & & & \Block{1-3}{\vspace*{8pt} n' \in \mathcal{N}_{v'}} & & \\
        \Block{3-1}{ n \in \mathcal{N}_{v} \hspace*{10pt}} & d_{n_1,n_1} & \dots & d_{n_1,n_j} & \infty & \dots & \infty \\
        & \vdots & \ddots & \vdots & \vdots & \ddots & \vdots \\
        & d_{n_i,n_1} & \dots & d_{n_i,n_j} & \infty & \dots & \infty \\ \hline
        \Block{3-1}{ n' \in \mathcal{N}_{v'} \hspace*{10pt}} & \infty & \dots & \infty & d_{n_{i'},n_{j'}} & \dots & d_{n_{i'},n_p} \\
        & \vdots & \ddots & \vdots & \vdots & \ddots & \vdots \\
        & \infty & \dots & \infty & d_{n_p,n_{i'}} & \dots & d_{n_p,n_p}
        \CodeAfter
            \OverBrace[shorten, yshift=5pt]{1-1}{1-3}{} 
            \OverBrace[shorten, yshift=5pt]{1-4}{1-6}{}
            \SubMatrix{\{}{1-1}{3-6}{.}[left-xshift=10pt]
            \SubMatrix{\{}{4-1}{6-6}{.}[left-xshift=10pt]
    \end{pNiceArray}.
\end{equation*}
Note that we ordered $\mathbf{D}$ as a block diagonal matrix here solely for better clarity, but this is not needed.

\section*{Acknowledgment}

The authors want to thank Fabian Neumann and the whole PyPSA team for their technical support during the model-building process and the integration of NPAP into PyPSA-Eur.

The authors used Claude (Anthropic) and Gemini (Google) to assist with writing the code for plots and tables, as well as for grammar checking and editorial refinement of the manuscript. The authors carefully reviewed the manuscript and take full responsibility for the content of this publication.

\printbibliography
    
\end{document}

